\documentclass{amsart}
\usepackage{graphicx}
\usepackage{amscd}
\usepackage{amsmath}
\usepackage{amsfonts}
\usepackage{amssymb}
\newtheorem{theorem}{Theorem}
\theoremstyle{plain}

\newtheorem{corollary}{Corollary}

\newtheorem{definition}{Definition}

\newtheorem{lemma}{Lemma}

\newtheorem{proposition}{Proposition}
\newtheorem{remark}{Remark}

\numberwithin{equation}{section}

\begin{document}
\title[On inclusions between ideals of polynomials]{On ideals of polynomials and their applications}
\author{Daniel M. Pellegrino}
\address[Daniel M. Pellegrino]{Depto de Matem\'{a}tica e Estat\'{\i}stica- Caixa Postal
10044- UFCG- Campina Grande-PB-Brazil }
\email{dmp@dme.ufcg.edu.br}
\thanks{This work is partially supported by Instituto do Mil\^{e}nio, IMPA}
\subjclass{Primary. 46G25, Secondary. 46B15}
\begin{abstract}In this paper we obtain some statements concerning ideals of polynomials and
apply these results in a number of different situations. Among other results,
we present new characterizations of $\mathcal{L}_{\infty}$-spaces, Coincidence
theorems, Dvoretzky-Rogers and Extrapolation type theorems for dominated polynomials.
\end{abstract}
\maketitle

%\subjclass{Primary.}
%
%

\section{Introduction, notations and background}

The notion of operators ideals goes back to Grothendieck \cite{Gro} and its
natural extension to polynomials and multilinear mappings is credited to
Pietsch \cite{Pietsch}. We prove some results on ideals of polynomials and
obtain, as corollaries and particular cases, new properties concerning
dominated, almost summing, integral polynomials and related ideals. Among
other results, we obtain Extrapolation and Dvoretzky-Rogers type theorems for
special ideals of polynomials, and prove new characterizations of
$\mathcal{L}_{\infty}$-spaces, extending results of Stegall-Retherford
\cite{Stegall} and Cilia-D'Anna-Guti\'{e}rrez \cite{Cilia}.

Throughout this paper $E,E_{1},...,E_{n},F,G,G_{1},...,G_{n},H$ will stand for
(real or complex) Banach spaces. Given \ a natural number $n\geq2,$ the Banach
spaces of all continuous $n$-linear mappings from $E_{1}\times...\times E_{n}$
into $F$ endowed with the $\sup$ norm will be denoted by $\mathcal{L}%
(E_{1},...,E_{n}$;$F)$ and the Banach space of all continuous $n$-homogeneous
polynomials $P$ from $E$ into $F$ with the $\sup$ norm is denoted by
$\mathcal{P}(^{n}E;F).$ If $T$ is a multilinear mapping and $P$ is the
polynomial generated by $T$, we write $P=\overset{\wedge}{T}.$ Conversely, for
the (unique) symmetric $n$-linear mapping associated to an $n$-homogeneous
polynomial $P$ we use the symbol $\overset{\vee}{P}.$ For $i=1,...,n,$
$\Psi_{i}^{(n)}:\mathcal{L}(E_{1},...,E_{n};F)\rightarrow\mathcal{L}%
(E_{i};\mathcal{L}(E_{1},\overset{[i]}{...},E_{n};F))$ will represent the
canonical isometric isomorphism defined by
\[
\Psi_{i}^{(n)}(T)(x_{i})(x_{1}\overset{[i]}{...}x_{n})=T(x_{1},...,x_{n}),
\]
where the notation $\overset{[i]}{...}$ means that the $i$-th coordinate is
not involved.

\begin{definition}
An ideal of multilinear mappings $\mathfrak{M}$ is a subclass of the class of
all continuous multilinear mappings between Banach spaces such that for all
index $n$ and $E_{1},$...,$E_{n},$ the components $\mathfrak{M}(E_{1}%
,...,E_{n};F)=\mathcal{L}(E_{1},...,E_{n};F)\cap\mathfrak{M}$ satisfy:

(i) $\mathfrak{M}(E_{1},...,E_{n};F)$ is a linear subspace of $\mathcal{L}%
(E_{1},...,E_{n};F)$ which contains the $n$-linear mappings of finite type.

(ii) If $A\in\mathfrak{M}(E_{1},...,E_{n};F),$ $u_{j}\in\mathcal{L}%
(G_{j};E_{j})$ for $j=1,...,n$ and $\varphi\in\mathcal{L}(F;H),$ then $\varphi
A(u_{1},...,u_{n})\in\mathfrak{M}(G_{1},...,G_{n};H).$

An ideal of (homogeneous) polynomials $\mathfrak{P}$ is a subclass of the
class of all continuous homogeneous polynomials between Banach spaces such
that for all index $n$ and all $E$ and $F$, \ the components $\mathfrak
{P}(^{n}E;F)=\mathcal{P}(^{n}E;F)\cap\mathfrak{P}$ satisfy:

(i) $\mathfrak{P}(^{n}E;F)$ is a linear subspace of $\mathcal{P}(^{n}E;F)$
which contains the polynomials of finite type.

(ii)If $P\in\mathfrak{P}(^{n}E;F),$ $\varphi_{1}\in\mathcal{L}(G;E)$ and
$\varphi_{2}\in\mathcal{L}(F;H),$ then $\varphi_{2}P\varphi_{1}\in\mathfrak
{P}(^{n}G;H).$
\end{definition}

There are several different ways to create ideals of polynomials and
multilinear mappings. We are mainly interested in two methods that we will
call (following \cite{B}) factorization and linearization:

\begin{itemize}
\item (Factorization method) If $\mathfrak{I}$ is an operator ideal, a given
$T\in\mathcal{L}(E_{1},...,E_{n};F)$ is of type $\mathcal{L}[\mathfrak{I}]$
($T\in\mathcal{L}_{\mathcal{L}[\mathfrak{I]}}(E_{1},...,E_{n};F)$) if there
are Banach spaces $G_{1},...,G_{n}$, linear operators $\varphi_{j}\in
\mathfrak{I}(E_{j};G_{j}),$ $j=1,...,n,$ and $R\in\mathcal{L}(G_{1}%
,...,G_{n};F)$ such that $T=R(\varphi_{1},...,\varphi_{n}).$ A given
$P\in\mathcal{P}(^{n}E;F)$ is of type $\mathcal{P}_{\mathcal{L}[\mathfrak{I}%
]}$ ($P\in\mathcal{P}_{\mathcal{L}[\mathfrak{I}]}(^{n}E;F)$) if there exist a
Banach space $G$, a linear operator $\varphi\in\mathfrak{I}(E;G)$ and
$Q\in\mathcal{P}(^{n}G;F)$ such that $P=Q\varphi.$
\end{itemize}

It is well known that $P\in\mathcal{P}_{\mathcal{L}(\mathfrak{I})}%
(^{n}E;F)\Leftrightarrow\overset{\vee}{P}\in\mathcal{L}_{\mathcal{L}%
(\mathfrak{I})}(^{n}E;F).$

\begin{itemize}
\item (Linearization method) If $\mathfrak{I}$ is an operator ideal, a given
$T\in\mathcal{L}(E_{1},...,E_{n};F)$ is of type $[\mathfrak{I}]$
($T\in\mathcal{L}_{[\mathfrak{I}]}(E_{1},...,E_{n};F)$) if $\Psi_{i}%
^{(n)}(T)\in\mathfrak{I}(E_{i};\mathcal{L}(E_{1},\overset{[i]}{...},E_{n}))$
for every $i=1,...,n.$ We say that $P\in\mathcal{P}(^{n}E;F)$ is of type
$\mathcal{P}_{[\mathfrak{I]}}$ ($P\in\mathcal{P}_{[\mathfrak{I}]}(^{n}E;F)$)
if $\overset{\vee}{P}$ is of type $[\mathfrak{I}]$.
\end{itemize}

The classes $\mathcal{P}_{\mathcal{L}[\mathfrak{I}]}$ and $\mathcal{P}%
_{[\mathfrak{I}]}$ are ideals of polynomials (see \cite{B},\cite{Braunss}%
,\cite{G}). For details on ideals of operators we refer to \cite{Pietsch2},
and for the theory of polynomials on infinite dimensions we mention
\cite{Dineen} and \cite{Mujica}.

If $\mathfrak{I}$ is the ideal $as,p$ of absolutely $p$-summing operators, by
a standard use of Ky Fan's Lemma it is not hard to obtain the following
characterization for the symmetric $n$-linear mappings of type $[as,p]$:

\begin{proposition}
\label{propn}A continuous (symmetric) multilinear mapping
$T:E\times...\times E\rightarrow F$ is of type $[as,p]$ if, and
only if,\ \ there exist $C\geq0$ and a regular probability measure
$\mu\in P\left(  B_{E^{\prime}}\right)$ such that
\begin{equation}
\left\|  T\left(  x_{1},...,x_{n}\right)  \right\|  \leq C\left\|
x_{1}\right\|  ...\left\|  x_{n-1}\right\|  \left[  \int_{B_{E^{\prime}}%
}\left|  \varphi\left(  x_{n}\right)  \right|  ^{p}d\mu\left(  \varphi\right)
\right]  ^{\frac{1}{p}}.
\end{equation}
\end{proposition}

\section{Preliminary results}

The well known Grothendieck-Pietsch Domination Theorem tells that a continuous
linear operator $T:E\rightarrow F$ is absolutely $p$-summing if, and only if,
there exist $C\geq0$ and a regular probability measure $\mu\in P\left(
B_{E^{\prime}}\right)  $ such that
\begin{equation}
\left\Vert T\left(  x\right)  \right\Vert \leq C\left[  \int_{B_{E^{\prime}}%
}\left\vert \varphi\left(  x\right)  \right\vert ^{r}d\mu\left(
\varphi\right)  \right]  ^{\frac{1}{r}} \label{Pit}%
\end{equation}
for every $x$ in $E$.

The concept of $p$-dominated polynomials (multilinear mappings) is one of the
most natural generalizations of absolutely summing operators and has been
broadly investigated (more information can be found in \cite{irish}%
,\cite{Matos2},\cite{Tonge},\cite{studia},\cite{P3}). We say that an
$n$-homogeneous polynomial is said to be $p$-dominated if there exist $C\geq0$
and a regular probability measure $\mu\in P\left(  B_{E^{\prime}}\right)  ,$
such that
\begin{equation}
\left\Vert P\left(  x\right)  \right\Vert \leq C\left[  \int_{B_{E^{\prime}}%
}\left\vert \varphi\left(  x\right)  \right\vert ^{p}d\mu\left(
\varphi\right)  \right]  ^{\frac{n}{p}}.
\end{equation}
We write $\mathcal{P}_{d,p}(^{n}E;F)$ to denote the space of $p$-dominated
$n$-homogeneous polynomials from $E$ into $F.$ It is well known that
$\mathcal{P}_{d,p}(^{n}E;F)=\mathcal{P}_{\mathcal{L}[as,p]}(^{n}E;F)$.

The following simple lemma, which proof we omit, will be useful later:

\begin{lemma}
\label{aaa}If $\mathfrak{I}_{1}$ and $\mathfrak{I}_{2}$ are ideals of
polynomials, and $\mathcal{L}_{\mathfrak{I}_{1}}(E;F)\subset\mathcal{L}%
_{\mathfrak{I}_{2}}(E;F)$ for every $F,$ then, for every $m$ and every $F$,%
\[
\mathcal{P}_{\mathcal{L}[\mathfrak{I}_{1}]}(^{m}E;F)\subset\mathcal{P}%
_{\mathcal{L}[\mathfrak{I}_{2}]}(^{m}E;F)\text{ and }\mathcal{P}%
_{[\mathfrak{I}_{1}]}(^{m}E;F)\subset\mathcal{P}_{[\mathfrak{I}_{2}]}%
(^{m}E;F).
\]
In particular, if $\mathcal{L}_{as,p}(E;F)=\mathcal{L}_{as,q}(E;F)$ for every
$F,$ then, for every $m$ and every $F$,%
\[
\mathcal{P}_{d,p}(^{m}E;F)=\mathcal{P}_{d,q}(^{m}E;F)\text{ and }%
\mathcal{P}_{[as,p]}(^{m}E;F)=\mathcal{P}_{[as,q]}(^{m}E;F)\text{ }.
\]
\end{lemma}

The next two propositions generalize \cite[Theorems 16 and 17]{Tonge}.

\begin{proposition}
If $E$ is a Banach space with cotype $2$ then for any Banach space $F,$ every
$n$ and every $p\leq2,$ we have
\begin{equation}
\mathcal{P}_{[as,p]}(^{n}E;F)=\mathcal{P}_{[as,2]}(^{n}E;F)\text{ and
}\mathcal{P}_{d,p}(^{n}E;F)=\mathcal{P}_{d,2}(^{n}E;F).\label{b3}%
\end{equation}
\end{proposition}

Proof. From a linear result of Maurey (see \cite[Theorem 36]{Woy}) we have
$\mathcal{L}_{as,p}(^{n}E;F)=\mathcal{L}_{as,2}(^{n}E;F)$ for every Banach
space $F$. Call on Lemma \ref{aaa} and obtain (\ref{b3}).$\Box$

\begin{proposition}
If $E$ is a Banach space with cotype $q,$ $2<q<\infty,$ then for any Banach
space $F$ and every $n,$ we have
\begin{equation}
\mathcal{P}_{[as,r]}(^{n}E;F)=\mathcal{P}_{[as,1]}(^{n}E;F)\text{ and
}\mathcal{P}_{d,r}(^{n}E;F)=\mathcal{P}_{d,1}(^{n}E;F)\label{b33}%
\end{equation}
for all $1<r<q^{\prime}$, where $\frac{1}{q}+\frac{1}{q^{\prime}\ }=1.$
\end{proposition}

Proof. A linear result due to Maurey (see \cite[Corollary 11.6]{Diestel})
asserts that if $E$ has cotype $q,$ $2<q<\infty,$ then
\[
\mathcal{L}_{as,r}(^{n}E;F)=\mathcal{L}_{as,1}(^{n}E;F)
\]
for all $1<r<q^{\prime}$ and every $F$. Again, we just need to use Lemma
\ref{aaa}.$\Box$

The following proposition will play an important role in several situations.
In particular, we will obtain new characterizations of $\mathcal{L}_{\infty}%
$-spaces, extending results of Stegall-Retherford \cite{Stegall} and
Cilia-D'Anna-Guti\'{e}rrez \cite{Cilia}, Extrapolations Theorems and results
of Dvoretzky-Rogers type for some classes of polynomials.

\begin{proposition}
\label{teoee}If $\mathcal{I}_{1}$ and $\mathcal{I}_{2}$ are ideals of
polynomials such that $\mathcal{P}_{\mathcal{I}_{1}}(^{n}E;F)\subset
\mathcal{P}_{\mathcal{I}_{2}}(^{n}E;F)$ for some Banach spaces $E$ and $F$,
some natural number $n$ and suppose that the following hold true:

(i)If $P\in\mathcal{P}_{\mathcal{I}_{2}}(^{n}E;F),$ then $\overset{\vee}%
{P}(.,a,...,a)\in$ $\mathcal{L}_{\mathcal{I}_{2}}(E;F)$ for every $a\in E$, fixed.

(ii) If $P\in\mathcal{P}_{\mathcal{I}_{1}}(^{m}E;F)$ and $\varphi
\in\mathcal{L}(E;\mathbb{K}),$ then $P.\varphi\in\mathcal{P}_{\mathcal{I}_{1}%
}(^{m+1}E;F)$, for all $m<n.$

Then $\mathcal{L}_{\mathcal{I}_{1}}(E;F)\subset\mathcal{L}_{\mathcal{I}_{2}}(E;F).$
\end{proposition}

Proof. If $T\in\mathcal{L}_{\mathcal{I}_{1}}(E;F)$, then define $\varphi\in$
$\mathcal{L}(E;\mathbb{K}),$ $\varphi\neq0$ and $a\in E$ such that
$\varphi(a)=1.$ Consider the following $n$-homogeneous polynomial:
\[
R(x)=\varphi(x)^{n-1}T(x).
\]
By applying (ii), $R\in\mathcal{P}_{\mathcal{I}_{1}}(^{n}E;F)\subset
\mathcal{P}_{\mathcal{I}_{2}}(^{n}E;F).$ Thus, (i) yields that $\overset{\vee
}{R}(.,a,...,a)\in\mathcal{L}_{\mathcal{I}_{2}}(E;F)$ and hence
\[
\frac{1}{n}T+\frac{n-1}{n}T(a)\varphi\in\mathcal{L}_{\mathcal{I}_{2}}(E;F).
\]
Since $\varphi\in\mathcal{L}(E;\mathbb{K})=$ $\mathcal{L}_{\mathcal{I}_{2}%
}(E;\mathbb{K}),$ we obtain $T\in\mathcal{L}_{\mathcal{I}_{2}}(E;F).\Box$

We have the following straightforward consequences:

\begin{corollary}
\label{ccc}If $\mathcal{I}_{1}$ and $\mathcal{I}_{2}$ are ideals of
polynomials such that $\mathcal{P}_{\mathcal{I}_{1}}(^{n}E;F)=\mathcal{P}%
_{\mathcal{I}_{2}}(^{n}E;F)$ for some Banach spaces $E$ and $F,$ some natural
number $n$ and suppose that the following hold true:

(i) If $P\in\mathcal{P}_{\mathcal{I}_{i}}(^{n}E;F),$ then $\overset{\vee}%
{P}(.,a,...,a)\in$ $\mathcal{L}_{\mathcal{I}_{i}}(E;F)$ for every $a\in E$,
fixed and $i=1,2.$

(ii) If $P\in\mathcal{P}_{\mathcal{I}_{i}}(^{m}E;F)$ and $\varphi
\in\mathcal{L}(E;\mathbb{K}),$ then $P.\varphi\in\mathcal{P}_{\mathcal{I}_{i}%
}(^{m+1}E;F)$, for all $m<n$ and $i=1,2.$

Then $\mathcal{L}_{\mathcal{I}_{1}}(E;F)=\mathcal{L}_{\mathcal{I}_{2}}(E;F).$
\end{corollary}

\begin{corollary}
\label{teoees}If $\mathcal{I}$ is an ideal of polynomials and $\mathcal{P}%
_{\mathcal{I}}(^{n}E;F)=\mathcal{P}(^{n}E;F)$ for some Banach spaces $E$ and
$F$, some natural number $n$ and
\[
P\in\mathcal{P}_{\mathcal{I}}(^{n}E;F)\Rightarrow\overset{\vee}{P}%
(.,a,...,a)\in\mathcal{L}_{\mathcal{I}}(E;F)
\]
for every $P\in\mathcal{P}_{\mathcal{I}}(^{n}E;F)$ and every $a\in E$, fixed,
then $\mathcal{L}_{\mathcal{I}}(E;F)=\mathcal{L}(E;F).$
\end{corollary}

Although Proposition \ref{teoee} and Corollaries \ref{ccc} and
\ref{teoees} suffice to our aims, it is worth remarking that it is
possible to obtain completer results, as we can see on the
following result, whose proof we omit.

\begin{theorem}
\label{t2}If $\mathcal{I}_{1}$ and $\mathcal{I}_{2}$ are ideals of polynomials
such that $\mathcal{P}_{\mathcal{I}_{1}}(^{n}E;F)\subset\mathcal{P}%
_{\mathcal{I}_{2}}(^{n}E;F)$ for some Banach spaces $E$ and $F$, some natural
number $n$ and suppose that the following hold true:

(i) If $P\in\mathcal{P}_{\mathcal{I}_{i}}(^{m}E;F),$ then $\overset{\vee}%
{P}(.,a^{k})\in$ $\mathcal{P}_{\mathcal{I}_{i}}(^{m-k}E;F)$ for all $m\leq n,$
every $k=1,...,m-1$, $a\in E$, fixed and $i=1,2$.

(ii) If $P\in\mathcal{P}_{\mathcal{I}_{i}}(^{m}E;F)$ and $\varphi
\in\mathcal{L}(E;\mathbb{K}),$ then $P.\varphi\in\mathcal{P}_{\mathcal{I}_{i}%
}(^{m+1}E;F)$, for all $m<n$ and $i=1,2.$

Then
\[
\mathcal{P}_{\mathcal{I}_{1}}(^{j}E;F)\subset\mathcal{P}_{\mathcal{I}_{2}%
}(^{j}E;F)
\]
for all $j=1,...,n.$ In particular, if $\mathcal{I}$ is an ideal of
polynomials and $\mathcal{P}_{\mathcal{I}}(^{n}E;F)=\mathcal{P}(^{n}E;F)$ for
some natural $n$ and

(i) If $P\in\mathcal{P}_{\mathcal{I}}(^{m}E;F),$ then $\overset{\vee}%
{P}(.,a^{k})\in\mathcal{L}_{\mathcal{I}}(E;F)$ for every $k=1,...,m-1,$ $a\in
E$, fixed, and $m\leq n.$

(ii) If $T\in\mathcal{L}_{\mathcal{I}}(^{m}E;F)$ and $\varphi\in
\mathcal{L}(E;\mathbb{K}),$ then $T.\varphi\in\mathcal{L}_{\mathcal{I}}%
(^{m+1}E;F)$, for all $m\leq n,$

then
\[
\mathcal{P}_{\mathcal{I}}(^{j}E;F)=\mathcal{P}(^{j}E;F)\text{ }\forall
j=1,...,n-1.
\]
\end{theorem}

Particular cases of this kind of structural property, called \ ``decreasing
scale property''\ (see G. Botelho \cite{B}), have been previously studied in
\cite{B}, \cite{Pellegrino}). In the next sections we will prove that several
ideals of polynomials satisfy the hypothesis of Theorem \ref{t2}. On the other
hand, it shall be emphasized that it is not difficult to find well known
ideals of polynomials which do not satisfy the decreasing scale property. In
fact, denoting by $(\mathcal{P}_{as,1}(^{n}E;F))_{n}$ the ideal of absolutely
$1$-summing polynomials endowed with the canonical norm (see \cite{Alencar} or
\cite{studia} for the definition of absolutely summing polynomials), one can
easily prove (exploring cotype properties) that $\mathcal{P}_{as,1}(^{n}%
l_{1};l_{1})=\mathcal{P}(^{n}l_{1};l_{1}),$for every $n\geq2$ and
$\mathcal{L}_{as,1}(l_{1};l_{1})\neq\mathcal{L}(l_{1};l_{1}).$

\section{Dvoretzky-Rogers and Extrapolation type Theorems}

An interesting consequence of Proposition \ref{teoee} is the following:

\begin{theorem}
(Dvoretzky-Rogers Theorem for $[al]$) If $E$ is a Banach space, then%
\[
\mathcal{P}_{[al]}(^{n}E;E)=\mathcal{P}(^{n}E;E)\Leftrightarrow\dim E<\infty.
\]
\end{theorem}

Proof. It is routine to verify that $[al]$ satisfy the hypothesis of
\ Corollary \ref{teoees}. Then we just need to call on the Dvoretzky-Rogers
Theorem for almost summing operators.$\Box$

Since, by Lemma \ref{aaa}, $\mathcal{P}_{d,p}(^{n}E;F)\subset\mathcal{P}%
_{[al]}(^{n}E;F),$ the Dvoretzky-Rogers Theorem for $[al]$ generalizes the
Dvoretzky Rogers Theorem for dominated polynomials, due to Matos
(\cite{Nachrichten}). For other concepts and results concerning almost summing
mappings we refer to \cite{JMAA}.

The next theorem lift to homogeneous polynomials a linear Extrapolation Theorem:

\begin{theorem}
(Extrapolation Theorems) Let $E$ be a Banach space and let $1\leq q.$ Suppose
that for \emph{some} natural $n$ and \emph{some} number $p$, $0<p<q$ we have%
\begin{equation}
\mathcal{P}_{d,p}(^{n}E;F)=\mathcal{P}_{d,q}(^{n}E;F)\text{ \emph{or}
}\mathcal{P}_{[as,p]}(^{n}E;F)=\mathcal{P}_{[as,q]}(^{n}E;F)\text{ for all
Banach spaces }F.\label{r1}%
\end{equation}
Then, for all Banach spaces $F,$\emph{ all} numbers $p$ such that $0<p<q$ and
\emph{all} natural $m,$ we have
\[
\mathcal{P}_{d,p}(^{m}E;F)=\mathcal{P}_{d,q}(^{m}E;F)\text{ \emph{and}
}\mathcal{P}_{[as,p]}(^{m}E;F)=\mathcal{P}_{[as,q]}(^{m}E;F).
\]
\end{theorem}

Proof. Firstly, we shall prove that ($\mathcal{I}_{1}=[as,p]$ , $\mathcal{I}%
_{2}=[as,q]$) and ($\mathcal{I}_{1}=d,p$ , $\mathcal{I}_{2}=d,q$) satisfy the
hypothesis of the Corollary \ref{ccc}. The proof that (ii) and (iii) of
Corollary \ref{ccc} hold for $d,p$ and $d,q$ is straightforward (it suffices
to call on the Grothendieck-Pietsch Domination Theorem). In order to verify
(i) for $[as,p]$ and $[as,q]$, we proceed as follows:

If $\ P\in\mathcal{P}_{[as,r]}(^{n}E;F)$ and $r=p$ or $q,$ then $\overset
{\vee}{P}\in\mathcal{L}_{[as,r]}(^{n}E;F)$ and thus there exist (from
Proposition \ref{propn}) $C\geq0$ and a regular probability measure $\mu\in
P\left(  B_{E^{\prime}}\right)  ,$ such that
\[
\left\|  \overset{\vee}{P}\left(  x_{1},...,x_{n}\right)  \right\|  \leq
C\left\|  x_{2}\right\|  ...\left\|  x_{n}\right\|  \left[  \int
_{B_{E^{\prime}}}\left|  \varphi\left(  x_{1}\right)  \right|  ^{r}d\mu\left(
\varphi\right)  \right]  ^{\frac{1}{r}}%
\]
and hence
\[
\left\|  \overset{\vee}{P}\left(  x,a,...,a\right)  \right\|  \leq C\left\|
a\right\|  ^{n-1}\left[  \int_{B_{E^{\prime}}}\left|  \varphi\left(  x\right)
\right|  ^{r}d\mu\left(  \varphi\right)  \right]  ^{\frac{1}{r}}.
\]
Thus $\overset{\vee}{P}\left(  .,a,...,a\right)  \in\mathcal{L}_{as,r}(E;F).$

For the proof of (ii) for $[as,r],$ and $r=p$ or $q,$ let us consider
$P\in\mathcal{P}_{[as,r]}(^{m}E;F)$ and $\varphi\in\mathcal{L}(E;\mathbb{K}).$
Then, defining the $(n+1)$-linear mapping $R(x_{1},...,x_{n})=\frac{1}%
{n+1}\varphi(x_{1})\overset{\vee}{P}(x_{2},...,x_{n+1})+...+$ $\frac{1}%
{n+1}\varphi(x_{n+1})\overset{\vee}{P}(x_{1},...,x_{n})$ we obtain%
\begin{align*}
&  \left(  \sum_{j=1}^{k}\left\|  \Psi_{1}^{(m+1)}(R)(x_{j})\right\|
^{r}\right)  ^{\frac{1}{r}}\\
&  \leq\frac{n}{n+1}\left(  \sum_{j=1}^{k}\left\|  \Psi_{1}^{(m)}%
(\overset{\vee}{P})(x_{j})\varphi\right\|  ^{r}\right)  ^{\frac{1}{r}}%
+\frac{1}{n+1}\left(  \sum_{j=1}^{k}\left\|  \varphi(x_{j})\overset{\vee}%
{P}\right\|  ^{r}\right)  ^{\frac{1}{r}}\\
&  \leq\frac{n}{n+1}\left\|  \varphi\right\|  \left(  \sum_{j=1}^{k}\left\|
\Psi_{1}^{(m)}(\overset{\vee}{P})(x_{j})\varphi\right\|  ^{r}\right)
^{\frac{1}{r}}+\frac{1}{n+1}\left\|  \overset{\vee}{P}\right\|  \left(
\sum_{j=1}^{k}\left\|  \varphi(x_{j})\right\|  ^{r}\right)  ^{\frac{1}{r}}\\
&  \leq\left(  \frac{n}{n+1}\left\|  \varphi\right\|  \left\|  \Psi_{1}%
^{(m)}(\overset{\vee}{P})\right\|  _{as,r}+\frac{1}{n+1}\left\|  \overset
{\vee}{P}\right\|  \left\|  \varphi\right\|  _{as,r}\right)  \left\|
(x_{j})_{j=1}^{k}\right\|  _{w,r}.
\end{align*}
Thus $R\in\mathcal{L}_{[as,r]}(^{m+1}E;F)$ and hence $\varphi P=\overset
{\wedge}{R}\in\mathcal{P}_{[as,r]}(^{m+1}E;F)$ and (ii) of Corollary \ref{ccc}
is satisfied.

Finally, from (\ref{r1}), by applying Corollary \ref{ccc} we obtain
\[
\mathcal{L}_{as,p}(E;F)=\mathcal{L}_{as,q}(E;F)
\]
for some $p$ such that $0<p<q$ and all Banach spaces $F.$ From a linear
Extrapolation Theorem due to Maurey (see \cite{Woy}), we conclude that
$\mathcal{L}_{as,p}(E;F)=\mathcal{L}_{as,q}(E;F)$ for all $p$ such that
$0<p<q$ and all Banach spaces $F.$ Call on Lemma \ref{aaa} to complete the
proof.$\Box$

A slightly different Extrapolation Theorem for dominated polynomials can be
found in \cite{Pellegrino} and \cite{P3}.

\section{Characterizations of $\mathcal{L}_{\infty}$-spaces}

The concept of $\mathcal{L}_{p}$-spaces, introduced by Lindenstrauss and Pe\l
czy\'{n}ski, in the seminal paper \textquotedblright Absolutely summing
operators in $\mathcal{L}_{p}$-spaces and their applications
\textquotedblright\ \cite{Pelcz} is a natural tool for a good understanding of
several properties of operators between Banach spaces. Since then,
characterizations of $\mathcal{L}_{p}$-spaces have been studied (see
\cite{Stegall},\cite{Cilia}, and \cite{Diestel} for other references). In this
section we will give new characterizations of $\mathcal{L}_{\infty}$-spaces.

Firstly let us recall the definition of integral polynomials, due to Cilia-D
'Anna-Guti\'{e}rrez \cite{Cilia}. An $m$-homogeneous polynomial
$P:E\rightarrow F$ is said to be integral if there exists $C\geq0$ such that
\[
\left|  \sum\limits_{i=1}^{n}<\varphi_{i},P(x_{i})>\right|  \leq C\sup
_{\psi\in B_{E^{\prime}}}\left\|  \sum\limits_{i=1}^{n}\left[  \psi
(x_{i})\right]  ^{m}\varphi_{i}\right\|
\]
for every natural $n,$ every $(x_{i})_{i=1}^{n}$ in $E$ and all $(\varphi
_{j})_{j=1}^{n}$ in the dual of $F$.

Analogously, an $n$-linear mapping $T\in\mathcal{L}(E_{1},...,E_{n};F)$ is
said to be integral if there exists $C\geq0$ such that%
\[
\left|  \sum\limits_{i=1}^{n}<\varphi_{i},T(x_{i}^{(1)},...,x_{i}%
^{(m)})>\right|  \leq C\underset{k=1,...,m}{\underset{\psi_{k}\in
B_{E_{k}^{\prime}}}{\sup}}\left\|  \sum\limits_{i=1}^{n}\left[  \psi_{1}%
(x_{i}^{(1)})...\psi_{m}(x_{i}^{(m)})\right]  \varphi_{i}\right\|
\]
for every natural $n,$ every $(x_{i})_{i=1}^{n}$ in $E$ and all $(\varphi
_{j})_{j=1}^{n}$ in the dual of $F$. It is not hard to see that the integral
polynomials forms an ideal of polynomials and from now on it will be denoted
by $\mathcal{I}\mathfrak{.}$ If $P$ is an integral ($n$-homogeneous)
polynomial from $E$ into $F$ we write $P\in\mathcal{P}_{\mathcal{I}}(^{n}E;F).$

\begin{lemma}
\label{coro} If $\mathcal{P}_{d,1}(^{m}E;F)\subset\mathcal{P}_{\mathcal{I}%
}(^{m}E;F),$ $\mathcal{P}_{[as,1]}(^{m}E;F)\subset\mathcal{P}_{\mathcal{I}%
}(^{m}E;F),$ $\mathcal{P}_{[as,1]}(^{m}E;F)\subset\mathcal{P}_{[\mathcal{I}%
]}(^{m}E;F)$ or $\mathcal{P}_{d,1}(^{m}E;F)\subset\mathcal{P}_{[\mathcal{I}%
]}(^{m}E;F),$ then
\[
\mathcal{L}_{as,1}(E;F)\subset\mathcal{L}_{\mathcal{I}}(E;F).
\]
\end{lemma}

Proof. The part $\mathcal{P}_{d,1}(^{m}E;F)\subset\mathcal{P}_{\mathcal{I}%
}(^{m}E;F)\Rightarrow\mathcal{L}_{as,1}(E;F)\subset\mathcal{L}_{\mathcal{I}%
}(E;F)$ is proved, using tensor products, in \cite{Cilia}. Here we give a
different and non-tensorial proof. We just need to verify that the ideals of
polynomials $d,1,$ $[as,1]$ satisfy (ii) and $[\mathcal{I}]$ and $\mathcal{I}$
satisfy (i) of Proposition \ref{teoee}. The verification for $d,1$ and
$[as,1]$ is already done in the proof of the Extrapolation Theorem. Concerning
the ideal of integral polynomials, an adequate handling of the polarization
formula yields to conclude that
\[
P\in\mathcal{P}_{\mathcal{I}}(^{m}E;F)\Leftrightarrow\overset{\vee}{P}%
\in\mathcal{L}_{\mathcal{I}}(^{m}E;F).
\]
Hence, if $P\in\mathcal{P}_{\mathcal{I}}(^{m}E;F)$, then
\[
\left|  \sum\limits_{i=1}^{n}<\varphi_{i},\overset{\vee}{P}(x_{i}%
^{(1)},...,x_{i}^{(m)})>\right|  \leq C\underset{k=1,...,m}{\underset{\psi
_{k}\in B_{E_{k}^{\prime}}}{\sup}}\left\|  \sum\limits_{i=1}^{n}\left[
\psi_{1}(x_{i}^{(1)})...\psi_{m}(x_{i}^{(m)})\right]  \varphi_{i}\right\|  .
\]
Choosing $(x_{i}^{(1)})_{i=1}^{n}=(a,...,a)$,..., $(x_{i}^{(m-1)})_{i=1}%
^{n}=(a,...,a),$ we have
\begin{align*}
\left|  \sum\limits_{i=1}^{n}<\varphi_{i},\overset{\vee}{P}(a,...,a,x_{i}%
^{(m)})>\right|   &  \leq C\underset{k=1,...,m}{\underset{\psi_{k}\in
B_{E_{k}^{\prime}}}{\sup}}\left\|  \sum\limits_{i=1}^{n}\left[  \psi
_{1}(a)...\psi_{m-1}(a)\psi_{m}(x_{i}^{(m)})\right]  \varphi_{i}\right\| \\
&  \leq C\left\|  a\right\|  ^{m-1}\sup_{\psi\in B_{E^{\prime}}}\left\|
\sum\limits_{i=1}^{n}\psi(x_{i}^{(m)})\varphi_{i}\right\|
\end{align*}
and thus $\overset{\vee}{P}(a,...,a,.)$ is an integral operator and
$\mathcal{P}_{\mathcal{I}}(^{m}E;F)$ satisfies (ii) of Proposition \ref{teoee}.

If $P\in\mathcal{P}_{[\mathcal{I}]}(^{m}E;F)$, then $\Psi_{1}^{(m)}%
(\overset{\vee}{P})\in\mathcal{L}_{\mathcal{I}}(E;\mathcal{L}(^{m-1}E;F)).$
Hence, if $(\varphi_{i})_{i=1}^{n}$ are in $F^{\prime},$ $(x_{i})_{i=1}^{n}$
are in $E$ and $a\in E$, define $\omega_{i}\in(\mathcal{L}(^{m-1}%
E;F))^{\prime}$ by $\omega_{i}(T)=\varphi_{i}(T(a,...,a)).$ We thus have
\begin{align*}
\left|  \sum\limits_{i=1}^{n}<\varphi_{i},\overset{\vee}{P}(a,...,a,x_{i}%
)>\right|   &  =\left|  \sum\limits_{i=1}^{n}<\omega_{i},\Psi_{1}%
^{(m)}(\overset{\vee}{P})(x_{i})>\right| \\
&  \leq C\underset{}{\underset{\psi\in B_{E^{\prime}}}{\sup}}\left\|
\sum\limits_{i=1}^{n}\psi(x_{i})\omega_{i}\right\| \\
&  =C\underset{}{\underset{\psi\in B_{E^{\prime}}}{\sup}}\sup_{\left\|
T\right\|  =1}\left\|  \sum\limits_{i=1}^{n}\psi(x_{i})\omega_{i}(T)\right\|
\\
&  =C\underset{}{\underset{\psi\in B_{E^{\prime}}}{\sup}}\sup_{\left\|
T\right\|  =1}\left\|  \sum\limits_{i=1}^{n}\psi(x_{i})\varphi_{i}%
(T(a,...,a))\right\| \\
&  \leq C\left\|  a\right\|  ^{n}\underset{}{\underset{\psi\in B_{E^{\prime}}%
}{\sup}}\left\|  \sum\limits_{i=1}^{n}\psi(x_{i})\varphi_{i}\right\|
\end{align*}
and the proof is done.$\Box$

A linear Theorem due to Stegall-Retherford \cite{Stegall} states that a Banach
space $E$ is an $\mathcal{L}_{\infty}$-space if and only if $\mathcal{L}%
_{as,1}(E;F)\subset\mathcal{L}_{\mathcal{I}}(E;F).$ In a recent paper,
Cilia-D'Anna-Guti\'{e}rrez \cite{Cilia} extend the characterization of
Stegall-Retherford and prove that a Banach space $E$ is an $\mathcal{L}%
_{\infty}$-space if and only if $\mathcal{P}_{as,1}(^{n}E;F)\subset
\mathcal{P}_{\mathcal{I}}(^{n}E;F)$ for some (every) natural $n.$ But, as it
can be seen in the next result, we can push a little further:

\begin{theorem}
Let $E$ be a Banach space. The following are equivalent:

(i) $E$ is an $\mathcal{L}_{\infty}$-space;

(ii) for all natural $n$ and every $F,$ we have $\mathcal{P}_{d,1}%
(^{n}E;F)\subset\mathcal{P}_{\mathcal{L}[\mathcal{I}\mathfrak{]}}(^{n}E;F);$

(iii) there is a natural $n$ such that for every $F,$ we have $\mathcal{P}%
_{d,1}(^{n}E;F)\subset\mathcal{P}_{\mathcal{L}[\mathcal{I}\mathfrak{]}}(^{n}E;F);$

(iv) for all natural $n$ and every $F,$ we have $\mathcal{P}_{d,1}%
(^{n}E;F)\subset\mathcal{P}_{[\mathcal{I}\mathfrak{]}}(^{n}E;F);$

(v) there is a natural $n$ such that for every $F,$ we have $\mathcal{P}%
_{d,1}(^{n}E;F)\subset\mathcal{P}_{[\mathcal{I}\mathfrak{]}}(^{n}E;F);$

(vi) for all natural $n$ and every $F$, $\mathcal{P}_{[as,1]}(^{n}%
E;F)\subset\mathcal{P}_{[\mathcal{I}]}(^{n}E;F);$

(vii) there is a natural $n$ such that for every $F,$ we have $\mathcal{P}%
_{[as,1]}(^{n}E;F)\subset\mathcal{P}_{[\mathcal{I}]}(^{n}E;F);$

(viii) for all natural $n$ and every $F,$ we have $\mathcal{P}_{d,1}%
(^{n}E;F)\subset\mathcal{P}_{\mathcal{I}}(^{n}E;F);$

(ix) there is a natural $n$ such that for every $F,$ we have $\mathcal{P}%
_{d,1}(^{n}E;F)\subset\mathcal{P}_{\mathcal{I}}(^{n}E;F);$
\end{theorem}

Proof. If $E$ is an $\mathcal{L}_{\infty}$-space, by the linear
characterization of $\mathcal{L}_{\infty}$-spaces, we have $\mathcal{L}%
_{as,1}(E;F)\subset\mathcal{L}_{\mathcal{I}}(E;F),$ for every Banach space
$F$. Now, since $\mathcal{P}_{d,1}(^{n}E;F)=\mathcal{P}_{\mathcal{L}%
[as,1]}(^{n}E;F)$, Lemma \ref{aaa} furnishes $\mathcal{P}_{d,1}(^{n}%
E;F)\subset\mathcal{P}_{\mathcal{L}[\mathcal{I}\mathfrak{]}}(^{n}E;F)$ for
every $n,$ and hence (i)$\Rightarrow$(ii) is done.

(ii)$\Rightarrow$(iii) is obvious. In order to prove (iii)$\Rightarrow$(i) we
just need to observe that we always have $\mathcal{P}_{\mathcal{L}%
[\mathcal{I}\mathfrak{]}}(^{n}E;F)\subset\mathcal{P}_{[\mathcal{I]}}(^{n}%
E;F)$. We thus have $\mathcal{P}_{d,1}(^{n}E;F)\subset\mathcal{P}%
_{[\mathcal{I]}}(^{n}E;F)$ and hence (by Lemma \ref{coro})%
\[
\mathcal{L}_{as,1}(E;F)\subset\mathcal{L}_{\mathcal{I}}(E;F),
\]
and consequently the linear characterization of $\mathcal{L}_{\infty}$-space
asserts that $E$ is an $\mathcal{L}_{\infty}$-space.

(ii)$\Rightarrow$(iv) is clear, since $\mathcal{P}_{\mathcal{L}[\mathcal{I}%
\mathfrak{]}}(^{n}E;F)\subset\mathcal{P}_{[\mathcal{I}]}(^{n}E;F).$
(iv)$\Rightarrow$(v) is obvious and (v)$\Rightarrow$(i) is a straightforward
consequence of Lemma \ref{coro}.

For the proof of (i)$\Rightarrow$(vi), if $E$ is an $\mathcal{L}_{\infty}%
$-space, then $\mathcal{L}_{as,1}(E;F)\subset\mathcal{L}_{\mathcal{I}}(E;F)$
for every $F$ and hence Proposition \ref{aaa} furnishes $\mathcal{P}%
_{[as,1]}(^{n}E;F)\subset\mathcal{P}_{[\mathcal{I}\mathfrak{]}}(^{n}E;F)$ for
every $n$.

(vi)$\Rightarrow$(vii) is obvious. For (vii)$\Rightarrow$(i), we call on Lemma
\ref{coro}, obtain $\mathcal{L}_{as,1}(E;F)\subset\mathcal{L}_{\mathcal{I}%
}(E;F)$ and one more time the linear characterization of $\mathcal{L}_{\infty
}$-space yields that $E$ is an $\mathcal{L}_{\infty}$-space.

Since $\mathcal{P}_{\mathcal{L}[\mathcal{I}\mathfrak{]}}(^{n}E;F)\subset
\mathcal{P}_{\mathcal{I}}(^{n}E;F)$ (\cite[Corollary 2.7]{Cilia}) it is clear
that (iii)$\Rightarrow$(viii). The proof of (viii)$\Rightarrow$(ix) is obvious
and we obtain (ix)$\Rightarrow$(i) by invoking Lemma \ref{coro}.

\begin{remark}
It is relevant to verify that, for example, $\mathcal{P}_{[as,1]}(^{n}E;F)$
and $\mathcal{P}_{d,1}(^{n}E;F)$ are different spaces, in general. In fact, it
can be proved (using a characterization of Hilbert Schmidt operators due to
Pe\l czy\'{n}ski \cite{HS}), following a suggestion of M. C. Matos, that
$P:l_{2}\rightarrow\mathbb{\ K}$ given by $P(x)=%
%TCIMACRO{\dsum \limits_{j=1}^{\infty}}%
%BeginExpansion
{\displaystyle\sum\limits_{j=1}^{\infty}}
%EndExpansion
\frac{1}{j^{\alpha}}x_{j}^{2}$ with $\alpha=\frac{1}{2}+\varepsilon$ and
$0<\varepsilon<\frac{1}{2}$ is such that $P\in\mathcal{P}_{[as,1]}(^{2}%
l_{2};\mathbb{K})$ and
$P\notin\mathcal{P}_{d,1}(^{2}l_{2};\mathbb{K})$. The proof is
given in \cite{Pell2} and \cite{Pellegrino}, but we sketch the
reasoning, for completeness.
\end{remark}

One can verify that $\overset{\vee}{P}:l_{2}\times
l_{2}\rightarrow\mathbb{K}\ $\ is given by
$\overset{\vee}{P}(x,y)=\sum\limits_{j=1}^{\infty}\frac{1}{j^{\alpha
}}x_{j}y_{j}$ and $(\frac{1}{j^{\alpha}})_{j=1}^{\infty}\in
l_{2}.$

It suffices to prove that $\overset{\vee}{P}$ fails to be
$1$-dominated, and
$\Psi_{1}^{(2)}(\overset{\vee}{P})\in\mathcal{L}_{as,1}(l_{2};l_{2}).$
Since
\[
\left(  \sum\limits_{j=1}^{m}\left\Vert \overset{\vee}{P}(e_{j},e_{j}%
)\right\Vert ^{\frac{1}{2}}\right)  ^{2}=\left[
\sum\limits_{j=1}^{m}\left( \frac{1}{j^{\alpha}}\right)
^{\frac{1}{2}}\right]  ^{2}\geq\left[ \sum\limits_{j=1}^{m}\left(
\frac{1}{m^{\frac{\alpha}{2}}}\right)  \right] ^{2}=m^{2-\alpha},
\]
if we had
\[
\left(  \sum\limits_{j=1}^{m}\left\Vert \overset{\vee}{P}(e_{j},e_{j}%
)\right\Vert ^{\frac{1}{2}}\right)  ^{2}\leq C\left\Vert (e_{j})_{j=1}%
^{m}\right\Vert _{w,1}^{2},
\]
we would obtain $m^{2-\alpha}\leq C(m^{\frac{1}{2}})^{2}=Cm $ and
it is a contradiction since $\alpha<1.$

In order to prove that $\Psi_{1}^{(2)}(\overset{\vee}{P})\in\mathcal{L}%
_{as,1}(l_{2};l_{2}),$ we shall note that
\[
\Psi_{1}^{(2)}(\overset{\vee}{P})((x_{j})_{j=1}^{\infty})=\left(
\frac {1}{j^{\alpha}}x_{j}\right)  _{j=1}^{\infty}.
\]
Now, a result of Pe\l czy\'{n}ski
(see \cite{Pelcz}) asserts that it suffices to show that $\Psi_{1}%
(\overset{\vee}{P})$ is a Hilbert-Schmidt operator. But is is easy
to check,
since $\sum\limits_{k=1}^{\infty}\left\Vert \Psi_{1}^{(2)}(\overset{\vee}%
{P})(e_{k})\right\Vert
_{l_{2}}^{2}=\sum\limits_{k=1}^{\infty}\left[  \frac
{1}{k^{\alpha}}\right]  ^{2}<\infty.$

\end{document}